\documentclass[10pt]{article}
\usepackage{amsmath}
\usepackage{amssymb}
\usepackage{amsthm}
\usepackage{graphicx}
\usepackage{verbatim}

\font\elevenss=cmss11

\font\eightss=cmss8

\font\sixss=cmss8 at 6pt

\newfam\ssfam
\textfont\ssfam=\elevenss \scriptfont\ssfam=\eightss
 \scriptscriptfont\ssfam=\sixss

\theoremstyle{plain}
\newtheorem{thm}{Theorem}[section]

\newtheorem{lem}[thm]{Lemma}
\newtheorem{pr}[thm]{Proposition}
\newtheorem{cor}[thm]{Corollary}
\newtheorem{defn}[thm]{Definition}

\newtheorem{conj}[thm]{Conjecture}

\newtheorem{example}[thm]{Example}
\theoremstyle{remark}

\newtheorem*{unremark}{Remark}

\newtheorem*{unremarks}{Remarks}

\newcommand{\Em}[1]{\textbf{#1}}

\def\Cox{\qed}

\def\Z{\mathbb{Z}}

\def\R{\mathbb{R}}

\def\CC{{\cal C}}

\def\ee{\epsilon}
\def\E{{\mathbb E}}
\def\P{{\mathbb P}}
\def\Cox{\hfill \Box}

\def\disp{\displaystyle}
\def\one{{\bf 1}}
\def\|{{\, | \, }}
\def\F{{\mathcal F}}
\def\B{{\mathcal B}}

\def\PZ{\P_{\uparrow}}

\def\yy{{\bf y}}

\def\ZZ{{\cal Z}}
\def\WW{{\cal W}}
\def\cL{{\cal L}}
\def\covers{\; \cdot \kern-5pt >}
\def\mcovers{\triangleright}

\def\mean{\mu}
\def\msr{\nu}
\def\msrr{\rho}

\makeatletter \def\romenumi{ \def\theenumi{\roman{enumi}}
\def\p@enumi{\theenumi} \def\labelenumi{(\@roman\c@enumi)}}
\makeatother

\begin{document}
\begin{titlepage}
\begin{center}
{\large \bf Concentration of Lipschitz functionals of determinantal
and other strong Rayleigh measures}
\\[5ex]
\end{center}

\begin{flushright}
Robin Pemantle\footnote{Department of Mathematics,
University of Pennsylvania,
209 South 33rd Street,
Philadelphia, PA 19104, USA, pemantle@math.upenn.edu}$^,$
\footnote{Research supported in part by NSF grant \# DMS-0905937}
and
Yuval Peres\footnote{Microsoft Research,
1 Microsoft Way, Redmond, WA, 98052, USA,
peres@microsoft.com}
\end{flushright}

\vfill

\noindent{\bf Abstract:}
Let $\{ X_1 , \ldots , X_n \}$ be a collection of binary valued
random variables and let $f : \{ 0 , 1 \}^n \to \R$ be a Lipschitz
function.  Under a negative dependence hypothesis known as the
{\em strong Rayleigh} condition, we show that $f - \E f$
satisfies a concentration inequality.
The class of strong Rayleigh measures includes determinantal
measures, weighted uniform matroids and exclusion measures; some
familiar examples from these classes are generalized negative binomials
and spanning tree measures.  For instance, any Lipschitz-1 function
of the edges of a uniform spanning tree on vertex set $V$ (e.g.,
the number of leaves) satisfies the Gaussian concentration inequality
$\disp{\P (f - \E f \geq a) \leq \exp \left ( - \frac{a^2}{8 \, |V|}
\right ) }$.  We also prove a continuous version for concentration of
Lipschitz functionals of a determinantal point process.
\vfill



\end{titlepage}

\setcounter{equation}{0}
\section{Introduction}
\label{sec:intro}

Our goal in this paper is to prove concentration inequalities
for Lipschitz functions of certain collections of negatively
dependent binary-valued random variables.  To illustrate our
general methods we state our main result in a special case
that was motivated by a question of E.\ Mossel (personal communication).
\begin{thm} \label{th:sp tree}
Let $G = (V,E)$ be a finite connected graph, let $\P$ be the
uniform measure on the spanning trees of $G$, and for $e \in E$
let $X_e$ be the indicator function of the event that $e$ is
in the chosen spanning tree.  Let $f : \{ 0 , 1 \}^E \to \R$
be any function with Lipschitz constant~1.  Then
$$\P (f - \E f \geq a) \leq \exp \left ( - \frac{a^2}{8 \, |V|}
   \right ) \, .$$
\end{thm}
For example we might take $f$ to be one half the number of
vertices whose degree in the random tree is odd.  This result
is a consequence of more general results stated in
Section~\ref{sec:results} and Section~\ref{sec:examples}.

\subsection{Classical concentration inequalities}

Let $\{ X_n : n \geq 1 \}$ be independent Bernoulli random variables
with respective means $\{ p_n \}$.  Let $S_n := \sum_{k=1}^n X_k$
denote the partial sums, $\mean_n := \E S_n = \sum_{k=1}^n p_k$ denote
the means and $V_n := \sum_{k=1}^n p_k (1-p_k)$ denote the variance
of $S_n$.  The simple and well known one-sided tail estimate for
$S_n$ is the classical Gaussian bound
\begin{equation} \label{eq:chernoff}
\P (S_n - \mean_n \geq a) \leq
   \exp \left  (- \frac{2 a^2}{n} \right ) \, .
\end{equation}
Replacing $X_n$ with $1-X_n$ gives the two-sided bound
\begin{equation} \label{eq:two-sided}
\P (|S_n - \mean_n| \geq a) \leq 2
   \exp \left  (- \frac{2 a^2}{n} \right ) \, .
\end{equation}
The bound~\eqref{eq:chernoff} may be found, among other places,
in~\cite[Corollary~5.2]{mcdiarmid89}.
The references given there include~\cite[(2.3)]{hoeffding}
as well as~\cite{chernoff}, which proves the result for
identically distributed variables.

When $p_n$ and $1 - p_n$ are bounded away from zero, the
variance of $S_n$ is of order $n$ and this kind of bound
is the best one can expect.  However, when $n \gg \mean_n$,
one might hope for uniformity in $n$ via bounds in which
the exponent depends on $\mean_n$ and not on $n$.  For example,
if $\max_{j \leq n} p_j$ is small then $S_n$ is well approximated
by a Poisson variable with mean $\mu_n$.
The upper tail of a Poisson is not as thin as a Gaussian,
being $\exp [ - \Theta (a \log (a/\mean))]$ rather than
$\exp [ - \Theta (a^2 / \mean)]$.  The bound
\begin{equation} \label{eq:pois}
\P (S_n \geq a + \mean) \leq
   e^a \; \left ( \frac{\mean}{a + \mean} \right )^{a + \mean} \leq
   \exp \left [ - \frac{a^2}{2 (a + \mean)} \right ]
\end{equation}
is proved in~\cite[Theorem~1]{hoeffding}
and asymptotically matches the Poissonian upper tail.

\subsection{Generalizations} \label{ss:gen}

Our aim is to generalize~\eqref{eq:chernoff} or its Poissonian
version~\eqref{eq:pois} in two ways.
Instead of $S_n$ we consider arbitrary Lipschitz functions of
$X_1 , \ldots , X_n$, and instead of independent Bernoullis
we consider a more general negatively dependent collection of
binary random variables.  We will give a number of applications,
but before this, we briefly discuss what is known about each
of the two generalizations separately.

For the first generalization, let $\B_n$ denote the rank-$n$
Boolean lattice $\{ 0 , 1 \}^n$ and let $f : \B_n \to \R$ be
Lipschitz with respect to the Hamming distance.  Replacing
$f$ by $f/c$ if necessary, we will lose no generality in
assuming our Lipschitz functions to have Lipschitz constant~1,
and we do so hereafter; thus $|f(x) - f(x')| \leq 1$ whenever
$x$ and $x'$ are two strings differing in only one position.
When $\P$ is a product measure, a well known generalization
of~\eqref{eq:chernoff}~\cite{mcdiarmid89} is
\begin{equation} \label{eq:gen 1}
\P (f - \E f \geq a) \leq e^{-2 a^2 / n} \, .
\end{equation}

For the second generalization, we say that a collection of
random variables $\{ X_j \}$ in $\{ 0 , 1 \}$ is
\Em{negative cylinder dependent} if
\begin{equation} \label{eq:cyl 1}
\P (X_j = 1 \mbox{ for all } j \in S) \leq \prod_{j \in S} p_j
\end{equation}
and
\begin{equation} \label{eq:cyl 2}
\P (X_j = 0 \mbox{ for all } j \in S) \leq \prod_{j \in S} (1 - p_j) \, .
\end{equation}
Negative cylinder dependence implies the
inequalities~\eqref{eq:chernoff}--\eqref{eq:two-sided} (see,
e.g.,~\cite[Theorem~3.4]{panconesi-srinivasan} with $\lambda = 1$).
Lyons~\cite[Section~6]{lyons03} lists extensions and applications
including one to balls in bins~\cite{dubhashi-ranjan} and one to
determinantal measures~\cite{soshnikov00}.

It is not known whether these two generalizations can be combined.
The random variables $\{ X_j \, : \, 1 \leq j \leq n \}$ are said to be
\Em{negatively associated} if $\E fg \leq (\E f) (\E g)$ for
every pair $f,g$ of increasing functions on $\{ 0 , 1 \}^n$ such that
$f (X_1 , \ldots , X_n)$ depends only on the values $\{ X_i : i \in S \}$
and $g (X_1 , \ldots , X_n)$ depends only on the values
$\{ X_i : i \notin S \}$, for some subset $S \subseteq \{ 1 ,
\ldots , n \}$.  By induction, this implies the weaker property
of negative cylinder dependence.  E. Mossel (personal communication, 2009)
asked us whether the following holds.
\begin{conj} \label{conj:both}
Let $X_1 , \ldots, X_n$ be negatively associated binary-valued
random variables.  Let $f : \{ 0,1 \}^n \to \R$ be Lipschitz-1
and denote $f_n := f(X_1 , \ldots , X_n)$.
Then~\eqref{eq:gen 1} holds with the bound $\exp (- 2 a^2 / n)$
replaced by $c_0 \exp (- c \, a^2 / n)$ for some positive
constants $c_0$ and $c$.
\end{conj}

To see why the exponent must be weakened consider the example
of Bernoulli random variables $X_1 , \ldots , X_n$ with
$n$ even, $\{ X_1 , \ldots , X_{n/2} \}$ independent with
mean $1/2$, and $X_{n/2 + j} = 1 - X_j$ for $1 \leq j \leq n/2$.
These are negatively associated and yet the Lipschitz-1 function
$f := \sum_{j=1}^{n/2} X_j - \sum_{j=n/2+1}^n X_j$ has tail
probabilities on the order of $e^{- a^2 / n}$.  It is
possible that this is the worst example and that the
conjecture holds with $c=1$, but a resolution of the conjecture
would be interesting even without the optimal value of $c$.
A recent paper~\cite{farc08} appears to settle this conjecture
and more, but the relevant result in that paper, Theorem~2,
is not correct, the proof therein failing at equation~(6).

Recent investigations of negative dependence properties
indicate that negative association may not be sufficiently
robust to use as a hypothesis in this context.  The problem
was posed in~\cite{pemantle-NA} to find a more useful
negative dependence property; this was answered
in~\cite{borcea-branden-liggett}, who showed that the
\Em{strong Rayleigh property} implies negative association
and many other desirable consequences, and is stable under
probabilistic operations such as conditioning, symmetrizing
and reweighting.

Our main result implies that Conjecture~\ref{conj:both} holds
with $c = 1/8$ if one assumes the strong Rayleigh property
rather than just negative association.
The strong Rayleigh property is known to hold for
most standard examples in which negative association is known
to hold, so this gives up little generality, and moreover the
strong Rayleigh property is usually easier to check than is
negative association.  Indeed for some of the measures described
below, the only way we know they are negatively associated is
by establishing the strong Rayleigh property.  Several classes
of measures satisfying the strong Rayleigh property are:
\begin{itemize}
\item Determinantal measures and point processes;
\item Bernoullis conditioned on the sum;
\item Measures obtained by running exclusion dynamics from a
deterministic starting state (or more generally, exclusion with
birth and death).
\end{itemize}

An overview of the rest of the paper is as follows.  In the
next section we introduce the strong Rayleigh property
and discuss its consequences.  One important consequence for
us will be the \Em{stochastic covering property}, which is
all we use to derive our basic concentration inequality.
In Section~\ref{sec:results} we state our results, and these
are proved in Section~\ref{sec:proofs}.  Section~\ref{sec:examples}
contains a number of applications.

\section{Strong Rayleigh property, stochastic covering property,
and other negative dependence conditions}
\label{sec:covering}

Let $[n]$ denote $\{ 1 , \ldots , n \}$ and let $\B_n := \{ 0 , 1 \}^n$
denote the Boolean lattice of rank $n$, with coordinatewise partial
order.  The function $N : \B_n \to \Z^+$ will be used throughout to
denote the counting function defined by $N(\omega) := \sum_{j=1}^n
\omega_j$.  A measure $\msr$ on $\B_n$ is said to be \Em{$k$-homogeneous}
if $\msr$ is supported on the set of $\{ \omega : N(\omega) = k \}$.
The probability
measure $\msr$ on $\B_n$ is said to be \Em{negatively associated}
if $\int fg \, d\msr \leq (\int f \, d\msr) (\int g \, d\msr)$ for
every pair of nonnegative monotone functions $f$ and $g$ such
that for some set $S \subseteq [n]$, the function $f$ depends only
on coordinates $\{ \omega_j : j \in S \}$ while the function $g$
depends only on coordinates $\{ \omega_j : j \in S^c \}$.

The \Em{strong Rayleigh condition} is said to hold for a measure $\P$
on $\B_n$ if the generating function
$$\sum_{\omega \in \B_n} \P (\omega) \prod_{j=1}^n z_j^{\omega_j}$$
has no roots $(z_1 , \ldots , z_n)$ all of whose coordinates lie
in the (strict) upper half plane.
This and many consequences are given in~\cite{borcea-branden-liggett},
including (implicitly) the \Em{stochastic covering property}
(see~Proposition~\ref{pr:SR -> SCP}), which was
conjectured~\cite[Conjecture~9]{pemantle-NA} to follow from
something a little weaker.  Some of the relevant implications
are summarized in Figure~\ref{fig:implications} below.

The definition of the stochastic covering property requires a few
preliminary definitions.  Recall that a measure $\msr$ on a partially
ordered set is said to \Em{stochastically dominate} a measure $\msrr$,
denoted $\msr \succeq \msrr$, if $\msr (A) \geq \msrr (A)$ for every upwardly
closed set $A$.  An equivalent condition is that there exists a
coupling, that is a measure $Q$ on $\B_n \times \B_n$ with
respective marginals $\msr$ and $\msrr$, supported on the set
$\{ (x,y) : x \geq y \}$.  If $\P$ is a measure on $\B_n$ making
the coordinate variables $\{ X_i : 1 \leq i \leq n \}$ negatively
associated, an immediate consequence of negative association is that
the conditional measure $(\P \| X_n = 0)$ on $\B_{n-1}$ stochastically
dominates the conditional measure $(\P \| X_n = 1)$.

\begin{figure}[ht]
\centering
\hspace{0.3in}
\includegraphics[scale=0.50]{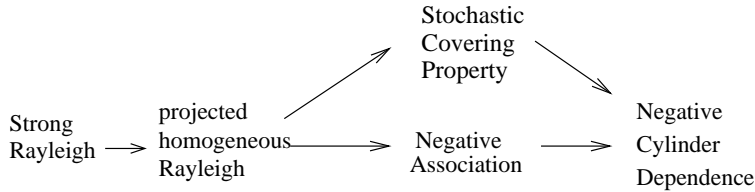}
\caption{relations among negative dependence properties}
\label{fig:implications}
\end{figure}

We say that the probability measure $\msr$ on $\B_n$ \Em{stochastically
covers} another probability measure $\msrr$ if there is a measure on
$\B_n^2$ with first
marginal $\msr$ and second marginal $\msrr$ (in other words, a coupling)
supported on the set of pairs $(x,y)$ for which $x = y$ or $x$ covers
$y$ in the coordinatewise partial order; here $x$ is said to cover $y$
when $x > y$ but there is no $z$ such that $x > z > y$.  We denote
the covering relation in $\B_n$ by $x \covers y$, and one measure
covering another by $\msr \mcovers \msrr$.  Stochastic covering is strictly
stronger than stochastic domination, and may be thought of as
``stochastic domination, but by at most~1''.

Stochastic covering
combines stochastic ordering with closeness in the so-called
$L^\infty$-transportation metric, defined on probability
measures on a given metric space as follows: $d_\infty
(\mu , \nu)$ is the least $\rho$ such that there is a coupling
of $\mu$ and $\nu$ supported on the set $\{ (x,y) : |x-y| \leq \rho \}$.
Thus $\mu \mcovers \nu$ implies $d_\infty (\mu , \nu) \leq 1$.  This
is useful because if $||f||_{\rm Lip}$ denotes the Lipschitz norm on
Lipschitz functions, then
\begin{equation} \label{eq:transport}
\left | \int f \, d\mu - \int f \, d\nu \right |
  \leq ||f||_{\rm Lip} \; d_\infty (\mu , \nu) \, .
\end{equation}

Suppose that $x \geq y$ and we compare the conditional laws $\P_x :=
(\P \| X_j = x_j, j \in S)$ and $\P_y := (\P \| X_j = y_j, j \in S)$
on the remaining coordinates, that is as laws on  $\{ 0 , 1 \}^{S^c}$.
If $\P$ and all its conditionalizations are negatively associated, it
follows that $\P_x \preceq \P_y$.

\begin{defn}[stochastic covering property]
We say that a probability measure $\msr$ on $\B_n$ has the
stochastic covering property if for every $S \subseteq \{ 1 ,
\ldots , n \}$ and for every $x,y \in \{ 0 , 1 \}^S$ with $x \covers y$,
the conditional law $(\msr \| X_j = x_j, j \in S)$ is covered by
the conditional law $(\msr \| X_j = y_j, j \in S)$.
\end{defn}

In~\cite[Theorem~4.2]{borcea-branden-liggett} it was shown that
strong Rayleigh property implies the \Em{projected homogeneous Rayleigh}
property (PHR), meaning that the measure can be embedded as the
first $n$ coordinates of a homogeneous measure $\msr'$ on
$\B_m$, for some $m \geq n$, that
has the ordinary Rayleigh property; the ordinary Rayleigh property
is that the partial derivatives of the generating function
$F(z_1 , \ldots , z_n) := \E \prod_{j=1}^n z_j^{X_j}$ satisfy
$F_i F_j \geq F_{ij} F$ at any point with positive real coordinates.
We record two further consequences.

\begin{pr} \label{pr:SR -> SCP}
PHR (and hence strong Rayleigh) implies the stochastic covering property.
\end{pr}

\noindent{\sc Proof:} PHR implies negative association of all
conditionalizations (CNA)~\cite[Theorem~4.10]{borcea-branden-liggett};
the homogeneous extension $\msr'$ witnessing the PHR property is
also PHR hence also CNA.  By negative association,
if $x \covers y$ then $(\msr' \| X_j = y_j , j \in S) \succeq
(\msr' \| X_j = x_j , j \in S)$, when viewed as measures on the
coordinates in $[m] \setminus S$.  Because $\msr'$ is homogeneous,
the coupling that witnesses this $\succeq$ relation in fact witnesses
the relation $\mcovers$.  Restricting to $[n] \setminus S$ we see
that $(\msr \| X_j = y_j , j \in S) \mcovers (\msr \| X_j = x_j , j \in S)$.
$\Cox$

\begin{pr}[\protect{\cite[Theorem~4.19]{borcea-branden-liggett}}]
\label{pr:k cover}
Let $\P$ be a strong Rayleigh measure on $\B_n$.
Let $\P_k$ denote $\P$ conditioned on $N = k$.  Then for every
$0 \leq k \leq n-1$ such that $\P (N=k)$ and $\P(N=k+1)$ are
both nonzero, we have the covering relation $\P_{k+1} \mcovers \P_k$.
$\Cox$
\end{pr}

\section{Results} \label{sec:results}

The chief consequence of the strong Rayleigh property that
we use to prove concentration inequalities is the
stochastic covering property.
Although all of our examples so far of measures with the SCP are
in fact strong Rayleigh, we note that this may not be the case in
the future, and with this in mind, we state a result that uses
only the SCP.
\begin{thm}[homogeneity and SCP implies Gaussian concentration]
\label{th:main}
Let $\P$ be a $k$-homogeneous probability measure on $\B_n$
satisfying the SCP.  Let $f$ be a Lipschitz-1 function on $\B_n$.
Then
$$ \P (f - \E f \geq a) \leq \exp \left ( - \frac{a^2}{8k} \right
   ) \, .  $$
\end{thm}
\begin{unremarks}
$(i)$ Replacing $f$ with $-f$ gives immediate two-sided bounds:
\begin{equation} \label{eq:two sided}
\P (|f - \E f| > a) \leq 2 \exp \left ( \frac{-a^2}{8 k} \right )  \, .
\end{equation}
$(ii)$ Replacing every $X_i$ by $1 - X_i$ if necessary, we may
assume without loss of generality that $k \leq n/2$, whence
\begin{equation} \label{eq:k to n}
\P (f - \E f > a) \leq \exp \left ( \frac{-a^2}{4 n} \right ) \, .
\end{equation}
\end{unremarks}

For strong Rayleigh measures that are not necessarily homogeneous,
we have the following result.
\begin{thm}[Gauss-Poisson bounds for general strong Rayleigh measures]
\label{th:general}
Let $\P$ be strong Rayleigh with mean $\mean = \E N$.
Let $f : \B_n \to \R$ be Lipschitz-1.  Then
\begin{eqnarray*}
\P (f - \E f > a) & \leq &
   3 \exp \left ( - \frac{a^2}{16 (a + 2 \mean)} \right ) \, ;  \\[1ex]
\P (|f - \E f| > a) & \leq &
   5 \exp \left ( - \frac{a^2}{16 (a + 2 \mean)} \right ) \, .
\end{eqnarray*}
\end{thm}
\begin{unremark}
Because $a, \mean \leq n$, the denominator in these inequalities
may be replaced by $48 n$.
\end{unremark}

\subsection*{Continuous versions}

Continuous versions of these results may be stated in terms of
point processes, which we now briefly review.  Formally, a point
process on a space $S$ is a random counting measure on $S$.
In other words, a point process is a map $\ZZ$ defined on a
probability space $(\Omega , \F , \P)$ taking values in the
space of counting measures on $S$, a counting measure being
one that takes only integer values or $+\infty$.
Intuitively, one envisions the sample counting measure
$\ZZ (\omega)$ as a set of points such that the sum of
delta functions at these points is the sample counting measure.

If the number $k$ of points in the support of $\ZZ$ is
deterministic, we may dispense with much of the formalism
by ordering the points in the support of $\ZZ$ uniformly
at random and identifying the process $\ZZ$ with the resulting
exchangeable probability law on sequences of length $k$ in $S$.
Notationally, if $\ZZ$ is a $k$-homogeneous point process on
$\R^d$ with law $\P$, we denote by $\PZ$ the corresponding
exchangeable law on $(\R^d)^k$.  For $1 \leq j \leq k$,
we use $X_j$ to denote the ``$j^{th}$ random point'', that
is, the $j^{th}$ coordinate function on $(\R^d)^k$.
The following sampling algorithm for any $k$-homogeneous
point process is almost trivial once one identifies $\ZZ$
with $\PZ$, and yet it is a generalization of an algorithm
previously proved only in the case of determinantal point
process in~\cite[Proposition~4.4.3]{HKPV}.

\begin{lem}[sampling in $k$ steps] \label{lem:k step}
Let $\ZZ$ be a $k$-homogeneous point process on a standard
Borel space $S$ and let $\PZ$ be the corresponding
exchangeable measure on $S^k$.  Then for $0 \leq j < k$
there are regular conditional distributions $Q_{x_1 , \ldots , x_j}$
for the law of $X_{j+1}$ given $X_1 = x_1 , \ldots , X_j = x_j$
such that the following procedure samples from $\PZ$.
\begin{quote}
Sample $X_1$ from $Q_\emptyset$.

Recursively, conditional on $X_1 = x_1 , \ldots , X_j = x_j$,
sample $X_{j+1}$ from $Q_{x_1 , \ldots , x_j}$.
\end{quote}
In the case where $S$ is finite, let $R$ denote the random
set $\{ X_1 , \ldots , X_n \}$.  Then the law $Q_{x_1 ,
\ldots , x_j}$ is equal to $1 / (k-j)$ times the conditional
intensity measure of $R \setminus \{ x_1 , \ldots , x_j \}$
given $x_1 , \ldots , x_j \in R$.
\end{lem}

\noindent{\sc Proof:}  Any standard Borel space admits
regular conditional distributions~\cite[Theorem~4.1.6]{durrett}.
The sampling algorithm essentially restates the definition
of regular conditional probabilities for sequential sampling.
Because $\PZ$ is exchangeable, conditioning on $X_{i_1} = x_1 ,
\ldots , X_{i_j} = x_j$ gives the same exchangeable measure
on the sequence of remaining elements of $R$ for any
$i_1 , \ldots , i_j$.  Thus for any $x$ other than
$x_1 , \ldots , x_j$, we have
$$\P (x \in R \| x_1 \in R , \ldots , x_j \in R) =
   (k-j) \P (X_{j+1} = x \| X_1 = x_1 , \ldots , X_j = x_j)$$
which is the final conclusion.
$\Cox$

\begin{unremark}
In the case of a measures on a finite set of size $n$, the
main point of this sampling scheme is to sample in $k$ steps
rather than $n$ steps, so as better to control the Azuma martingale.
But also, sequential conditioning on $x \in R$ can be easy to
compute.  For example, conditioning on an edge being in a
spanning tree replaces the original graph by a contraction
along that edge.
\end{unremark}

Intuitively, there is a stochastic covering property for
point processes defined to hold when conditioning on
the presence of a point depresses the process everywhere
else but by at most one point.  To make this definition
precise begin by extending the notion of one measure
stochastically covering another to point processes.
We say the point process $\ZZ$ stochastically covers
the point process $\WW$ if there is a coupling of these
two laws on counting measures supported on pairs $(\mu , \nu)$
such that $\mu = \nu$ or $\mu = \nu + \delta_x$ for some $x$.
Metrizing the space of finite counting measures on $S$ by
the total variation distance, we see as before that if
$\ZZ \mcovers \WW$ then $d_\infty (\cL (\ZZ) , \cL (\WW)) \leq 1$.

Next, given a $k$-homogeneous
point process $\ZZ$ on a space $S$, we let $\ZZ_{x_1 , \ldots , x_j}$
denote the $(k-j)$-homogeneous point process whose law
is the law of $\{ X_{j+1} , \ldots , X_k \}$ when sampling
according to the procedure in Lemma~\ref{lem:k step}
conditional on $X_1 = x_1 , \ldots , X_j = x_j$.
The $k$-homogeneous point process $\ZZ$ is said to have the
\Em{stochastic covering property} if
$$\ZZ_{x_1 , \ldots , x_j} \mcovers \ZZ_{x_1 , \ldots , x_{j+1}}$$
for all choices of $x_1 , \ldots , x_{j+1}$.  Note
that the left-hand side is $(k-j)$-homogeneous while
the right-hand side is $(k-j-1)$-homogeneous.
\begin{thm} \label{th:cont homog}
Let $\ZZ$ be a $k$-homogeneous point process on a standard
Borel space $S$ and
let $f$ be a Lipschitz-1 function (with respect to the total
variation distance) on counting measures with
total mass $k$ on $\R^d$.  If $\ZZ$ has the SCP, then
$$\P (f - \E f \geq a)
   \leq \exp \left ( - \frac{a^2}{8k} \right ) \, .$$
\end{thm}

For point processes that are not homogeneous, as in the
discrete case, we require more than the SCP.  Rather than
defining a notion of strong Rayleigh here, we will stick
to the case of determinantal point processes, this being
where all of our examples arise; see Section~\ref{s:continuous}
for definitions.
\begin{thm} \label{th:continuous}
Let $\ZZ$ be a determinantal point process with $\E N = \mean < \infty$.
Let $f$ be a Lipschitz-1 function on finite counting measures.  Then
\begin{eqnarray*}
\P (f - \E f \geq a )
   & \leq & 3 \exp \left ( - \frac{a^2}{16 (a + 2\mean)} \right )
   \, ; \\[1ex]
\P (|f - \E f \geq a )
   & \leq & 5 \exp \left ( - \frac{a^2}{16 (a + 2\mean)} \right ) \, .
\end{eqnarray*}
\end{thm}

\section{Proofs} \label{sec:proofs}

\subsection{The classical proofs}

To prove bounds such as~\eqref{eq:chernoff}, one obtains
an upper bound for $\E e^{\lambda S_n}$, and then applies
Markov's inequality, choosing $\lambda$ optimally.  Underlying
the bounds on $\E e^{\lambda S_n}$ are corresponding bounds
for compensated increments.  Let $\Delta$ denote a variable
with mean zero.  Three classical exponential bounds are as
follows.
\begin{eqnarray}
|\Delta| \leq 1 & \Rightarrow &
   \E e^{\lambda \Delta} \leq e^{\lambda^2 / 2}
   \label{eq:delta 1} \\[1ex]
\Delta \in [r,s] & \Rightarrow &
   \E e^{\lambda \Delta} \leq e^{\lambda^2 (s-r)^2 / 8}
   \label{eq:delta 2} \\[1ex]
\Delta \in [r,s] & \Rightarrow &
\E e^{\lambda \Delta}
   \leq \exp \left [ (e^\lambda - 1 - \lambda) \; |rs| \right ] \, ,
   \label{eq:delta 3}
\end{eqnarray}
when $|r-s| \leq 1$.
These are used together with the following two special cases
of Markov's inequality.
\begin{eqnarray}
\E e^{\lambda X} \leq e^{c \lambda^2 / 2} & \Longrightarrow
   & \P (X \geq a) \leq e^{-a^2 / (2c)}
   \label{eq:AH 2} \\[1ex]
\E e^{\lambda X} \leq e^{b \, (e^\lambda - \lambda - 1)}
   & \Longrightarrow &
   \P (X > a) \leq e^a \left ( \frac{b}{a+b} \right )^{a+b}
   \leq \exp \left [ - \frac{a^2}{2 (a+b)} \right ]
   \label{eq:pois2}
\end{eqnarray}

These inequalities have appeared many times in the literature.
Inequalities~\eqref{eq:delta 1} and~\eqref{eq:AH 2} constitute
the classical Azuma-Hoeffding inequality and imply
\begin{eqnarray}
\E e^{\lambda (S_n - \mean_n)} & \leq & e^{\lambda^2 n / 2}
   \, ; \label{eq:AH 1} \\[1ex]
\P (S_n - \mean_n \geq a) & \leq & e^{-a^2 / (2n)} \, .  \label{eq:AH 3}
\end{eqnarray}
This is valid for any martingale with differences bounded by~1;
an exposition can be found in~\cite[Theorem~7.2.1]{alon-spencer}.
The improvement to~\eqref{eq:delta 2} is present already
in~\cite{hoeffding}, though the exposition in~\cite{mcdiarmid89}
is clearer (see Lemma~5.8 therein).  When the increments of
$S_n - \mean_n$ are compensated Bernoullis, one may take $b-a = 1$
rather than~2, resulting in an improvement by a factor of~4
in the exponent,
\begin{equation} \label{eq:quadratic}
\E e^{\lambda (S_n - \mean_n)} \leq e^{\lambda^2 n / 8} \, .
\end{equation}
which together with~\eqref{eq:AH 2} yields~\eqref{eq:chernoff}.
Finally,~\eqref{eq:delta 3} and induction yield
\begin{equation} \label{eq:exp pois}
\E e^{\lambda (S_n - \mean_n)} \leq e^{(e^\lambda - \lambda - 1) \, V_n}
   \leq e^{(e^\lambda - \lambda - 1) \, \mean_n}
\end{equation}
where $V_n := \sum_{k=1}^n p_k (1 - p_k)$ is the variance of $S_n$;
together with~\eqref{eq:pois2} this implies~\eqref{eq:pois};
these results appear, for instance, in~\cite[(1.3)--(1.6)]{freedman75}.

To prove the generalization to Lipschitz functions, let
$$M_k := \E \left ( f(X_1 , \ldots , X_n) \| X_1 , \ldots , X_k
   \right ) - \E f(X_1 , \ldots , X_n) \, .$$
It is immediate that $\{ M_k \}$ is a martingale and that
conditional on $X_1 , \ldots , X_{k-1}$, the two possible values
of $M_k$ differ by at most~1.  Hence, conditional on
$X_1 , \ldots , X_{k-1}$, the increment $\Delta_k := M_k - M_{k-1}$
is constrained to an interval of length at most~1.
Applying~\eqref{eq:delta 2} then yields~\eqref{eq:gen 1}.

The extension of inequalities~\eqref{eq:chernoff}--\eqref{eq:two-sided}
to negatively cylinder dependent random variables is established
by examining the power series for $e^{\lambda S_n}$.  This
may be expanded into positive sums of expectations of products
of powers of the variables $\{ X_j : 1 \leq j \leq n \}$.
Negative cylinder dependence implies that
these are bounded from above by the corresponding products of
expectations.  Therefore,~\eqref{eq:quadratic} and~\eqref{eq:exp pois}
hold when the assumption of independence is replaced by negative cylinder
dependence, whence the probability inequalities~\eqref{eq:chernoff}
and~\eqref{eq:pois} hold as well.  This and more is shown
in~\cite[Theorem~3.4]{panconesi-srinivasan}, specializing
their more general negative cylinder property to $\lambda = 1$.
We remark that only the first inequality~\eqref{eq:cyl 1}
in the definition of negative cylinder dependence
is used to obtain bounds on $\E e^{\lambda S_n}$ for $\lambda > 0$,
which suffices for the upper tail bounds.  Lower tail bounds require
these inequalities for $\lambda < 0$, for which the second
inequality~\eqref{eq:cyl 2} is required.

\subsection{Proof of Theorems~\protect{\ref{th:main}}
and~\protect{\ref{th:cont homog}}}

Theorem~\ref{th:main} is a special case of Theorem~\ref{th:cont homog}.
This is because any probability measure $\mu$ on $\B_n$ may be viewed
as the law of a point process on the $n$ element set $[n]$, where
the random counting measure $\ZZ (\omega)$ is defined by
$\ZZ (\omega) (S) = \sum_{j \in S} \omega_j$.  Informally,
the points of the process are the coordinates of the ones
in the sample $\omega$.  With this interpretation, the SCP
on $\B_n$ is inherited from the SCP for the point process $\ZZ$,
whence Theorem~\ref{th:cont homog} with $S = \R$ (or any other
standard Borel space containing $[n]$) implies Theorem~\ref{th:main}.
It remains to prove Theorem~\ref{th:cont homog}.

Let $(\Omega , \F , \P)$ be a probability space on which is
constructed the generalized sampling scheme described in
Lemma~\ref{lem:k step}.  Let $\F_j := \sigma (X_1 , \ldots , X_j)$
and let
\begin{equation} \label{eq:mart}
M_j := \E (f \| \F_j) - \E f \, .
\end{equation}
denote the martingale of sequential revelation.  Applying
the method of bounded differences is now mostly a matter of
bookkeeping.  At a sample point where $X_i = x_i, 1 \leq i \leq k$,
the quantity $M_j$ may be written as the integral of $f$ against the
law of the point process
$$\ZZ_{x_1 , \ldots , x_j} + \sum_{i=1}^j \delta_{x_i} \, .$$
By the SCP, we have $\ZZ_{x_1 , \ldots , x_j} \mcovers
\ZZ_{x_1 , \ldots , x_{j+1}}$ whence
$$d_\infty \left ( \ZZ_{x_1 , \ldots , x_j} + \sum_{i=1}^j \delta_{x_i}
   \;\; , \;\; \ZZ_{x_1 , \ldots , x_{j+1}} + \sum_{i=1}^{j+1} \delta_{x_i}
   \right ) \leq 2 \, .$$
By the Lipschitz assumption on $f$, it follows that
$|M_{j+1} - M_j| \leq 2$.
We now apply the basic Azuma-Hoeffding inequality~\eqref{eq:AH 3}
to $\{ M_j / 2 \}_{1 \leq j \leq k}$ yielding
$$\P (f - \E f \geq a) = \P \left ( \frac{M_k}{2} > \frac{a}{2} \right )
   \leq \exp \left ( - \frac{a^2}{8 k} \right ) \, .$$
$\Cox$

\subsection{Proof of Theorem~\protect{\ref{th:general}}}

In this section we assume $\P$ is the law of a strong
Rayleigh measure on $\B_n$ with finite mean $\E N = \mean$.  
We also let $f$ denote
an arbitrary but fixed Lipschitz-1 function on configurations
and define a function $\phi$ on $\Z^+$ by
$$\phi (k) := \E (f \| N = k) \, .$$

\begin{lem} \label{lem:bernoulli}
The variable $N$ is distributed as the sum of independent Bernoullis.
\end{lem}

\noindent{\sc Proof:} In the definition of the strong Rayleigh
property, setting the variables $z_1 , \ldots , z_n$ equal
produces a univariate polynomial with no roots in the upper half
plane.  As pointed out at the beginning of Section~3
of~\cite{borcea-branden-liggett}, such a polynomial with
real coefficients must have all its roots real.  Since
the coefficients are nonnegative, this implies it is the
generating function for a convolution of Bernoullis.
$\Cox$

\begin{lem} \label{lem:N concentration}
The variable $N$ satisfies
$$\E e^{\lambda (N - \mean)} \leq \exp \left [ \mean (e^\lambda - 1 - \lambda)
   \right ]$$
and consequently for any $a > 0$,
$$\P \left ( N \geq \mean + \frac{a}{2} \right ) \leq
   \exp \left ( - \frac{(a/2)^2}{2 (\mean + a/2)} \right ) \, .$$
\end{lem}

\noindent{\sc Proof:} By Lemma~\ref{lem:bernoulli} $N$ is
distributed as the sum of independent Bernoullis, which
implies the first inequality; this implies the second inequality
by~\eqref{eq:pois2}.
$\Cox$

\begin{lem} \label{lem:phi lipschitz}
The function $\phi$ is Lipschitz-1.
\end{lem}

\noindent{\sc Proof:} By Proposition~\ref{pr:k cover} in the case
of strong Rayleigh measures on $\B_n$
 we know that
$\P_{k+1} \mcovers \P_k$.  By definition of
the stochastic covering relation, $(\phi (k+1) , \phi (k))$
may be written as $\E (f(\eta) , f(\xi))$ where $d(\eta , \xi) = 1$
almost surely.  The conclusion then follows from the fact that
$f$ is Lipschitz-1.
$\Cox$

\begin{lem} \label{lem:phi concentration}
The random variable $\phi (N)$ satisfies the concentration inequality
$$\E e^{\lambda (\phi (N) - \E \phi (N))} \leq
   e^{\mean (e^\lambda - 1 - \lambda)} \, .$$
Consequently, the upper tails of $\phi (N)$ obey the bound
$$\P (\phi (N) - \E \phi (N) > t) \leq
   e^{-\frac{t^2}{2 (t + \mean)}} \, .$$
\end{lem}

\noindent{\sc Proof:}  Pursuant to Lemma~\ref{lem:bernoulli},
let $\{ Y_j \}$ be a finite or countably infinite collection of
independent Bernoulli variables whose sum has the same law as $N$;
we may therefore prove the statements with $N$ replaced by
$\sum_j Y_j$.  Write $\phi (\sum_j Y_j) - \E \phi (\sum_j Y_j)$
as the final term of a martingale $\{ M_\ell \}$ where
$M_\ell := \E ( \phi (\sum_j Y_j) \| \F_\ell) - \E \phi$ and
$\F_\ell := \sigma (Y_1 , \ldots , Y_\ell)$.  If the number
of Bernoullis is infinite, the final term is a limit
almost surely and in $L^2$.  The martingale $\{ M_\ell \}$
is a binary martingale, meaning that conditional on
$\F_\ell$, the distribution of $M_{\ell+1}$ is concentrated on
two values.  In other words,
$$(M_{\ell+1} \| \F_\ell) = p \delta_r + (1-p) \delta_s$$
where $p$ is the mean of the Bernoulli variable $Y_{\ell+1}$.
More importantly, $r = \int f \, d\mu$ and $s = \int f \, d\nu$
where $\mu$ is the conditional law of $\sum_j Y_j$ given
the values of $Y_1 , \ldots , Y_\ell$ (which are measurable
with respect to $\F_\ell$) and given $Y_{\ell + 1} = 1$,
and $\nu$ is the conditional law of $\sum_j Y_j$ given
the values of $Y_1 , \ldots , Y_\ell$ and given $Y_{\ell + 1} = 0$.
Clearly $\mu$ and $\nu$ are probability measures on $\Z^+$
satisfying $d_\infty (\mu , \nu) \leq 1$, whence because $\phi$
is Lipschitz-1, we see that $|r-s| \leq 1$.  From~\eqref{eq:delta 3}
we then obtain
$$\E \left ( e^{\lambda (M_{\ell+1} - M_\ell)} \| \F_\ell \right )
   \leq \exp \left ( p (1-p) (e^\lambda - 1 - \lambda) \right ) \, .$$
The lemma follows by induction.
$\Cox$

\noindent{\sc Proof of Theorem~\protect{\ref{th:general}}}
The event $\{ f - \E f > a \}$ is contained in the union of three events:
$$\left \{ N > \mean + \frac{a}{2} \right \} \cup
   \left \{ \phi (N) - \E f > \frac{a}{2} \right \} \cup
   \left \{ f - \phi (N) > \frac{a}{2} \, , \, N \leq \mean +
      \frac{a}{2} \right \}  \, .$$
Thus $\P (f - \E f > a)$ is bounded above by the sum of the
corresponding probabilities.  Each of the first two pieces is
bounded above by $\exp [- a^2 / (4(a + 2\mean))]$: the first
follows from Lemma~\ref{lem:N concentration} and the second
uses Lemma~\ref{lem:phi concentration}, noting that $\E \phi
= \E \; \E (f \| N) = \E f$.  The last piece is bounded above
by $\exp [- a^2 / (16 (a + 2\mean))]$.  To see this, observe
that the measures $\P_k$ are all strong Rayleigh (this
is~\cite[Corollary~4.18]{borcea-branden-liggett}).
For any $k \leq \mean + a/2$, we can apply Theorem~\ref{th:main} to the
homogeneous measure $\P_k$, obtaining
$$\P \left ( f - \phi (N) > \frac{a}{2} \| N = k \right )
   \leq \exp \left ( - \frac{(a/2)^2}{8 k} \right )
   \leq \exp \left ( - \frac{a^2}{16 (a + 2 \mean)} \right ) \, .$$
Reassembling these gives the upper bound
$$\P \left ( f - \phi (N) > \frac{a}{2} \, , \, N \leq \mean +
   \frac{a}{2} \right )
   \leq \exp \left ( -\frac{a^2}{16 (a + 2 \mean)} \right ) \, .$$
This last piece has the worst bound; using it for all three
pieces gives the first inequality of the theorem; we remark that
the better upper bound of $2 \exp [- a^2 / (4(a + 2\mean))] +
\exp [- a^2 / (16 (a + 2\mean))]$ is in fact valid.

For the two-sided bound, we need to consider two more events
in addition to the three already considered, namely the events
$\{ \phi (N) - \E \phi (N) < -a/2 \}$ and $\{ f - \phi (N) < -a/2, N \leq
\mean + a/2 \}$.  The arguments for these two extra events are
exactly analogous to two of the three arguments we have already
seen, leading to a bound of $\exp [ - a^2 / (16 (a + 2 \mean)) ]$
for each of the two new summands and establishing the two-sided bounds.
$\Cox$

\clearpage
\section{Applications} \label{sec:examples}
\setcounter{equation}{0}

In this section we discuss some classes of measures known to satisfy
the hypotheses of our concentration results.  The following Venn
diagram gives a sense of how these classes intersect each other.
\begin{figure}[ht]
\centering
\hspace{0.3in}
\includegraphics[scale=0.50]{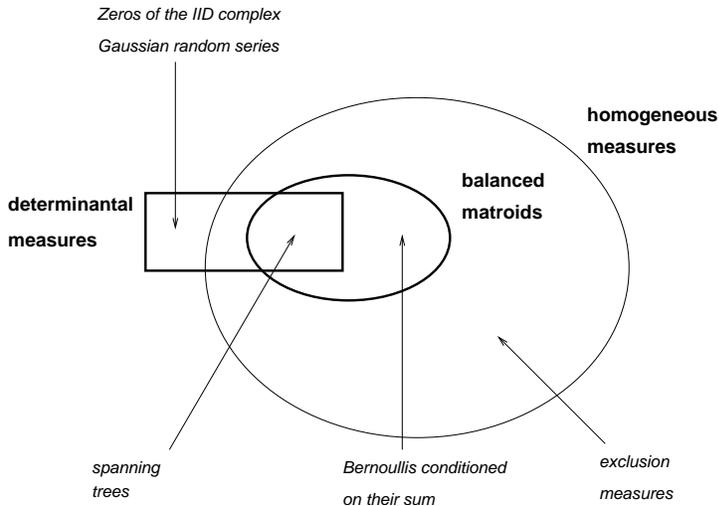}
\caption{some classes of strong Rayleigh measures}
\label{fig:venn}
\end{figure}

\subsection{Matroids} \label{ss:matroids}

A collection $\CC$ of subsets of a finite set $E$, all of a given
cardinality, $k$, is said to be the set of bases of a matroid
if it satisfies the base exchange axiom (see, e.g.,~\cite{tutte}):
if $A$ and $B$ are distinct members of $\CC$ and $a \in A
\setminus B$, then there exists $b \in B \setminus A$
such that $A \cup \{ b \} \setminus \{ a \} \in \CC$.
Given a matroid, it is natural to consider the uniform measure
on $\CC$.  More generally, the \Em{weighted random base} is
chosen from the probability measure
$$\msr_w (B) := C \prod_{e \in B} w(e) \, ,$$
where $\{ w(e) : e \in E \}$ is a collection of nonnegative
real numbers (weights) and $C$ is a normalizing constant.
Identifying $E$ with the set $\{ 1 , \ldots , |E| \}$, the
measure $\msr_w$ and the random variables $X_e := 1_{e \in B}$
can be thought of as living on $\B_{|E|}$.

For general matroids, $\E X_e X_f$ may be greater than
$(\E X_e) (\E X_f)$.  Some speculation has been given to the
most natural class of matroids for which negative correlation
or negative association must hold.  Feder and
Mihail~\cite{feder-mihail} define a balanced matroid
to be a matroid all of whose minors satisfy pairwise negative
correlation.  Their proof of the following fact was the
basis for the original proof of negative association for
determinantal processes~\cite[Theorem~6.5]{lyons03}.
\begin{pr}[\protect{\cite[Theorem~3.2]{feder-mihail}}] \label{pr:matroids}
The law $\msr_w$ of a random base of a balanced matroid,
multiplicatively weighted by the weighting function $w$,
has the SCP.
$\Cox$
\end{pr}

Because measures supported on the bases of a matroid are homogeneous,
there is nothing gained by improving the SCP to the strong Rayleigh
property, and we have the following immediate corollary.
\begin{cor} \label{cor:balanced}
Let $f$ be a Lipschitz-1 function with respect to Hamming
distance on the bases of a balanced matroid of rank $k$
on $n$ elements.  Then
$$\P (f - \E f > a) \leq \exp \left ( - \frac{a^2}{8k} \right ) \, .$$
$\Cox$
\end{cor}

\begin{example}[spanning trees] \label{eg:spanning tree}
One of the most important examples of a matroid is the
set of spanning trees of a finite, connected, undirected graph.
To spell this out, a spanning tree for a finite graph $G = (V,E)$
is a subset $E' \subseteq E$ such that $(V,E')$ is a connected
and acyclic.  The set of spanning trees is a matroid on $E$.
The weighted random spanning tree was shown to be a balanced
matroid by~\cite[Theorem~1]{burton-pemantle}.  In fact they
showed it is determinantal (see also~\cite[Example~1.1]{lyons03}
and~\cite[Example~4.3.2]{HKPV}), though at the time
consequences of being determinantal, such as
the Strong Rayleigh property, had not been developed.
Spanning trees are the only well known class of matroid
whose uniform (or weighted) measure is determinantal.

Let $f_0: \{ 0 , 1 \}^E \to \Z$ count the number of vertices
of odd degree in the graph defined by any subset of the edges.
Deleting or adding an edge changes $f_0$ by at most~2.
Let $f$ be the random variable resulting from applying
$(1/2) f_0$ to a the weighted random spanning tree on a
graph $G$.  Thus $f$ is a Lipschitz~1 function that counts
half the number of vertices that have odd degree in the
random tree.  Random variables that count local
properties such as this are of natural graph theoretic interest.
Parity counting variables similar to $f$ play a role, for example
in the randomized TSP approximation algorithm of~\cite{saberi}.
The number of edges in any spanning tree is $|V| - 1$.  An
application of Corollary~\ref{cor:balanced} immediately gives
the concentration inequality in Theorem~\ref{th:sp tree};
note that $|V|$, rather than $|E|$, appears in the denominator
of the exponent.
\end{example}

\begin{example}[conditioned Bernoullis, weighted matroids]
\label{eg:cond bern}
Let $\lambda_1 , \ldots , \lambda_n > 0$ be real numbers and let
$\mu$ be the measure on the subsets of $[n]$ of cardinality $k$
given by
$$\mu (x_1 , \ldots , x_n) = \frac{\prod_{j=1}^n \lambda_j^{x_j}}
   {\sum_{N (\yy) = k} \prod_{j=1}^n \lambda_j^{y_j}} \, .$$
We may think of $\mu$ in two ways.  The first is as a special
case of the weighted random base, specialized to $M(n,k)$,
the matroid whose bases are all the subsets of $[n]$
that have cardinality $k$.  The second is that it is the law of
independent Bernoulli variables $X_j$ with $\E X_j =
\lambda_j / (1 + \lambda_j)$, conditioned on $\sum_{j=1}^n X_j = k$.
We see from Proposition~\ref{pr:matroids} that $\mu$ is strong
Rayleigh.  Alternatively, we may deduce this from the strong
Rayleigh property for product measures along with closure under
conditioning on the sum~\cite[Corollary~4.18]{borcea-branden-liggett}.
Restricting to $[m]$, for $m < n$, gives a joint distribution
on $\B_m$ which may be thought of as a multivariate generalization
of the hypergeometric distribution.  Because the strong Rayleigh
property is inherited, this restriction is strong Rayleigh as well.
The resulting concentration properties of these measures
have been exploited in~\cite{saberi} in connection with
TSP approximation.  We remark that more general conditioning,
such as conditioning $N := \sum_j X_j$ to lie in an interval
of more than two points, does not preserve the strong Rayleigh
property.
\end{example}

\subsection{Exclusion measures} \label{ss:exclusion}

The symmetric group ${\cal S}_n$ acts on $\B_n$ by permuting
the coordinates.  Suppose a nonnegative rate $r(\tau)$
is given for each transposition $\tau \in {\cal S}_n$.
Define a random evolution on $\B_n$ by letting each
pair of coordinates $(i,j)$ transpose independently
at rate $r(\tau_{ij})$.  In other words, we have a
continuous time chain on $\B_n$ which jumps from $x$ to
$\tau (x)$ at rate $r(\tau)$ for each transposition $\tau$.
This process is known as the symmetric exclusion process.

Borcea, Branden and Liggett~\cite[Proposition~5.1]{borcea-branden-liggett}
prove that the strong Rayleigh property is preserved under this
evolution.  In particular, because the point mass at a single
state is always strong Rayleigh, it follows that the time $t$
distribution of a symmetric exclusion process started from a
deterministic state is strong Rayleigh.  The stochastic covering
property follows, as do PHR and negative association.  Interestingly,
before the publication of~\cite{borcea-branden-liggett}, all
that was known about this model was negative cylinder
dependence~\cite[Lemma~2.3.4]{liggett77}).

Recently, it was shown by~\cite{wagner-BAMS} that one can add
birth and death to the exclusion dynamics and still preserve
the strong Rayleigh property.  More specifically, let $\{ \alpha_i ,
\beta_i : 1 \leq i \leq n \}$ be positive real numbers and let
$\omega_i$ change to one at rate $\alpha_i$ and to zero at rate
$\beta_i$, along with the exclusion dynamics.  Then the evolution
preserves the strong Rayleigh property and in particular, if the
starting state is deterministic, all time $t$ marginals are
strong Rayleigh.

\begin{cor} \label{cor:exclusion}
Let $\P$ be the law on $\B_n$ resulting from running
an exclusion process for a fixed time, starting from
a deterministic state with $k$ sites occupied.  Then
$$\P (f - \E f \geq a) \leq e^{- a^2 / (8 k)} \, .$$
\end{cor}

\begin{example} \label{eg:exclusion}
Let $n > 0$ be even and populate a $n \times n$ square of the
integer lattice in $\Z^2$ (with torus boundary conditions)
by filling all sites in the left half and leaving empty
all sites in the right half.  Run the symmetric exclusion
process for time $t$ with rate~1 on each edge.  Let
$f_t (\omega)$ denote the number of edges at time
$t$ with exactly one endpoint occupied.  The mean of $f_t$
starts at $n$ at time 0 and approaches its limiting value
of $n^2 - O(1)$ as $t \to \infty$.  Once $t = \Theta (n^2)$,
the variance of $f_t$ becomes $\Theta (n^2)$ and the
concentration inequality
$$\P (f - \E f \geq a) \leq e^{- a^2 / (4 n^2 )} \, ,$$
which holds for all $t$, becomes a meaningful
Gaussian tail bound (here $k = n^2 / 2$).
\end{example}

\subsection{Determinantal measures on a finite Boolean lattice}
\label{ss:finite}

We say that a probability measure $\P$ on $\B_n$ is
\Em{determinantal} (in the general sense) if there is
an $n \times n$ real or complex matrix $K$ such that that
for every $S \subseteq \{ 1 , \ldots , n \}$,
\begin{equation} \label{eq:det}
\E \prod_{j \in S} X_j = \det K_S
\end{equation}
where $K_S$ is the submatrix of $K$ obtained by choosing only those
rows and columns whose index is in $S$.  In this definition, the
phrase ``general sense'' refers to the lack of further assumptions
on $K$.  An important subclass is the \Em{Hermitian} determinantal
measures, for which the matrix $K$ is Hermitian.  In this paper we
will be interested only in the Hermitian case and will use the term
\Em{determinantal} hereafter to refer only to the case where $K$ is
Hermitian.  Determinantal measures are known to be negatively
associated~\cite[Theorem~6.5]{lyons03}.  In fact they are strong
Rayleigh~\cite[proof of Theorem~3.4]{borcea-branden-liggett} and
therefore satisfy the stochastic covering property.
\begin{example}[uniform or weighted spanning tree]
\label{eg:spanning tree 2}
As previously remarked, the uniform or weighted random spanning
tree is a determinantal measure.
\end{example}

In the next section we will extend the notion of a determinantal
measure to the continuous setting.  The extension to a countably
infinite set of variables is more straightforward: the kernel $K$
is now indexed by a countably infinite set, but~\eqref{eq:det}
may be interpreted as holding for all finite sets $S$.  The following
example of a determinantal process on $\Z$ appeared first
in~\cite{johansson03}.

\begin{example}[positions of non-colliding RW's]
\label{eg:nonintersecting RW}
Let $\{ Y^{(k)} : 1 \leq k \leq n \}$ be $n$ independent time
homogeneous nearest neighbor random walks on $\Z$.  Begin the
walks at locations $y_1 , \ldots , y_n$ and suppose the
event that the walks are all at their starting positions at
time $2n$ and have not intersected has positive probability.
Conditional on this event, the positions at time $n$ form
a determinantal measure.  That is, the indicator functions
$\{ X_j \}$ have a determinantal law, where $X_j = 1$ if
some $Y^{(k)}$ is at position $j$ at time $n$, and zero otherwise.
\end{example}

\begin{unremark}
The positions of non-colliding random walks are given by a determinant
under more general conditions (see~\cite{karlin-mcgregor}).
The present situation is arranged so as
to make the kernel Hermitian.
\end{unremark}

\section{Determinantal point processes}
\label{s:continuous}

We consider here only simple point processes and often assume
$\E N < \infty$ as well.
If $\rho_k : (\R^d)^k \to \R^+$ are measurable functions, then the
simple point process $\ZZ$ is said to have joint intensities
$\{ \rho_k \}$ if for any $k$ and any family $D_1 , \ldots , D_k$
of disjoint Borel subsets of $\R^d$,
$$\E \left [ \prod_{j=1}^k \ZZ (D_j) \right ] = \int_{\prod_j D_j}
   \rho_k (x_1 , \ldots , x_k) \, dx_1 \cdots dx_k \, .$$
In particular,
$$\E N = \int_{\R^d} \rho_1 (x) \, dx$$
so under the assumption $\E N < \infty$, we see that
$\rho_1 (x) \, dx$ is a finite measure on $\R^d$.  If
$\rho_1$ is not finite, we will assume it is $\sigma$-finite.
In any case, $\rho_1$ is called the \Em{first intensity measure};
see~\cite[Sections~1.2~and~4.2]{HKPV} for further discussion of
joint intensities and determinantal measures.

\begin{defn}[determinantal point process]
A point process $\ZZ$ is said to be determinantal if it has
joint intensities $\{ \rho_k \}$ and there is a measurable kernel
$K : (\R^d)^2 \to \CC$ such that
\begin{equation} \label{eq:kernel}
\rho_k (x_1 , \ldots , x_k)
   = \det \left ( K(x_i , x_j) \right )_{1 \leq i,j \leq k} \, .
\end{equation}
If $K(y,x) = \overline{K(x,y)}$ for every $x,y$, then the process
is said to be Hermitian. When discussing determinantal processes below, we will always assume they are Hermitian.
\end{defn}

Stochastic covering carries over to the continuous case.
To state the relevant results we invoke the notion of
the \Em{Palm process}.  This is a version of the process
conditioned on the (measure zero) event of a point at a
specified location, $x$.  It may be obtained by conditioning
on there being a point within distance $\ee$ of a given
location $x$, then taking a weak limit.  A more complete
treatment may be found in~\cite{kallenberg}.  The   following proposition is proved
in~\cite{goldman10}.
\begin{pr}[\protect{\cite{goldman10}}] \label{pr:cover cont}
Suppose $\ZZ$ is a determinantal point process with continuous
kernel $K$ and finite trace.  Fix $x$ and let $\ZZ_x$ denote the
Palm process that conditions on a point at $x$.  Let $\ZZ_x'$
denote the result of removing the point at $x$ from $\ZZ_x$.
Then
\begin{enumerate} \romenumi
\item Whenever $K-L$ is positive semi-definite, the process with
kernel $K$ stochastically dominates the process with kernel $L$
(this is~\cite[Theorem~3]{goldman10}).
\item $\ZZ_x'$ is determinantal with kernel $L$ such that $K-L$ is
positive semi-definite.
\item Consequently, $\ZZ \succeq \ZZ_x'$.
\end{enumerate}~~\\[-4ex]
\end{pr}

   The continuous analogue
of Proposition~\ref{pr:k cover} is

\begin{pr} \label{eq:pr k cont}
Let $\ZZ$ be a determinantal point process with finite mean
$\E N = \mean < \infty$.  Then for any $k$ for which $\P (N = k+1)$
and $\P (N = k)$ are both nonzero, the conditional distributions
of $\ZZ$ given $N$ satisfy
$$( \ZZ \| N = k+1) \mcovers (\ZZ \| N = k) \, .$$
\end{pr}

\noindent{\sc Proof:}
The following facts may be found in~\cite[Theorem~7]{HKPV-PS}.
A determinantal point process $\ZZ$ with mean $\mean < \infty$
has a kernel $K$ whose spectrum is countable, contained in $[0,1]$,
and sums to $\mean$.  Furthermore, $\ZZ$ may be represented as a
mixture of homogeneous determinantal processes as follows.
Let $\{ \lambda_i : i \geq 1 \}$ enumerate the eigenvalues
with multiplicities and let $\{ \phi_i \}$ be a corresponding
eigenbasis.  For each $i$, flip an independent coin with success
probability $\lambda_i$.  Let $I$ denote the set of $i$ for which
the coin-flip was successful.  Let $K_I$ be the (random) projection
operator onto the subspace spanned by the eigenvectors $\phi_i$
for which the coin-flip was successful.  Then $K_I$ is almost
surely a projection of finite dimension $|I|$ and is the kernel
of a $|I|$-homogeneous determinantal point process.  Choosing
$K_I$ at random and then sampling from the corresponding process
recovers the law of $\ZZ$.

Several consequences are apparent.  First, conditioning on
$N = k$ is the same as conditioning on exactly $k$ successes
among the Bernoulli trials.  Secondly, the conditional law of $I$
given $|I| = k+1$ stochastically dominates the conditional
law of $I$ given $|I| = k$.  When the number of Bernoullis
is finite, this follows from the strong Rayleigh property
for independent Bernoullis; an easy limit argument extends
the conclusion to the infinite case.  This fact about stochastic
domination is equivalent to saying
that the conditional law of the random subspace $K_I$ given
$|I| = k+1$ stochastically dominates the conditional law
of the random subspace $K_I$ given $|I| = k$, in the
sense that the two laws can be coupled as $(K , K')$
so that $K' \subseteq K$.  When $K' \subseteq K$,
the operator $\pi_{K} - \pi_{K'}$ is positive semi-definite.
By~$(ii)$ of Proposition~\ref{pr:cover cont}, we conclude that
$(\ZZ \| N = k+1) \succeq (\ZZ \| N = k)$ which is equivalent
to stochastic covering in this case.
$\Cox$

\noindent{\sc Proof of Theorems~\protect{\ref{th:continuous}}:}

With $f$ as in the statement of the theorem, and $I$ the collection of indices described in the previous proposition, define
$\psi(I)$ to be the expectation of $f$ applied to a configuration chosen from the determinantal process with kernel $K_I$. Recall the notation $N=|I|$.

The event $\{ f - \E f > a \}$ is contained in the union of three events:
$$\left \{ N > \mean + \frac{a}{2} \right \} \cup
   \left \{ \psi(I) - \E f > \frac{a}{2} \right \} \cup
   \left \{ f - \psi(I) > \frac{a}{2} \, , \, N \leq \mean +
      \frac{a}{2} \right \}  \, .$$
Thus $\P (f - \E f > a)$ is bounded above by the sum of the
corresponding probabilities.  Each of the first two pieces is
bounded above by $\exp [- a^2 / (4(a + 2\mean))]$: the first
follows from Lemma~\ref{lem:N concentration} and the second
follows from the proof of Lemma~\ref{lem:phi concentration} because $\psi$ is Lipschitz in the Bernoulli variables $Y_i:=\one_{i\in I}$.
  The last piece is bounded above
by $\exp [- a^2 / (16 (a + 2\mean))]$.  To see this, apply Theorem~\ref{th:cont homog} to the homogenous determinantal processes $\P_I$ with kernels $K_I$, obtaining, when $|I|=k$, that
$$\P_I \left ( f - \psi(I) > \frac{a}{2} \right )
   \leq \exp \left ( - \frac{(a/2)^2}{8 k} \right )
   \leq \exp \left ( - \frac{a^2}{16 (a + 2 \mean)} \right ) \, .$$
Reassembling these gives the upper bound
$$\P \left ( f - \psi(I) > \frac{a}{2} \, , \, N \leq \mean +
   \frac{a}{2} \right )
   \leq \exp \left ( -\frac{a^2}{16 (a + 2 \mean)} \right ) \, .$$
The rest of the argument is identical to the conclusion of the proof of Theorem~\ref{th:general}.
$\Cox$

\begin{example}[Ginibre's translation invariant process]
\label{eg:ginibre}
Ginibre~\cite{ginibre65} considers the distribution of eigenvalues
of an $k \times k$ matrix with independent complex Gaussian
entries.  In the limit as $k \to \infty$, the density becomes
constant over the whole plane.  The limiting process $\ZZ$ turns
out to be a (Hermitian) determinantal point process with kernel
$$K(z_1 , z_2) := \frac{1}{\pi} e^{z_1 \overline{z_2}}
   \, \exp \left ( - \frac{|z_1|^2 + |z_2|^2}{2} \right ) \, ;$$
see, e.g.~\cite[(2.16)]{soshnikov00}.  The process $\ZZ$ is
ergodic and invariant under all rigid transformations of the
plane.  It was suggested~\cite{lecaer-ho90} to use this
process as the set of centers for a random Voronoi tesselation
because the mutual repulsion of the points makes the resulting
tesselation more realistic than the standard Poisson-Voronoi
tesselation for many purposes.  Some rigorous results along
these lines were obtained in~\cite{goldman10}.

The mean number of points in any region $D$ is $1/\pi$ times
the area $|D|$, so the restriction of $\ZZ_D$ to such
a region of finite area is a determinantal process with
finite mean number of points.
Fix a finite region, $D$, and let $f$ count the number of ``lonely''
points in $D$, these being such that no other point of $\ZZ$
in $D$ is within distance~1.  We claim that $f$ is Lipschitz
with constant equal to~6.  Clearly if a point $z$ is added
to the configuration $\eta$ then $f$ can increase by
at most~1.  It is well known that the maximum number
of points in a unit disk that can be at mutual distance
of at least~1 from one another is~6, which implies that
the addition of $z$ can result in the loss of at most~6
lonely points.  Applying Theorem~\ref{th:continuous}
to $f/6$ yields the concentration inequality
$$\P (|f - \E f| \geq a) \leq 5 \exp \left (
   - \frac{a^2}{96 (a + 12 |D| / \pi )} \right ) \, .$$
\end{example}

\begin{example}[Zeros of random polynomials] \label{eg:zeros}
Let $\{ X_n \}$ be IID standard complex Gaussian random variables
and define the random power series
$$h(z) := \sum_{n=0}^\infty X_n z^n \, .$$
It is easy to see that $h$ is almost surely analytic on the
open unit disk and the number of zeros on any disk of radius
$\rho < 1$ has finite mean.  The remarkable properties of the point
process $\ZZ$ on the unit disk that is the zero set of $h$ are
detailed in~\cite{peres-virag}.  It is a determinantal process
whose kernel is the Bergman kernel $\pi^{-1} (1 - z \overline{w})^{-2}$.
It is invariant under M\"obius transformations of the unit disk
and has intensity measure $\pi^{-1} / (1 - |z|^2)^2$.  Endowing
the unit disk with the hyperbolic metric, the M\"obius transformations
become isometries, whence $\ZZ$ is hyperbolic isometry invariant.

Fix $\rho < 1$ and $r > 0$ and let $f$ count the number of
zeros of the restriction $\ZZ_\rho$ of $\ZZ$ to the disk
of radius $1-\rho$ that are ``hyperbolically lonely'',
meaning that no other point of $\ZZ_\rho$ is within a
hyperbolic distance $r$.  Let $c_r$ denote the maximum
number of points at mutual hyperbolic distance $r$ that
may be be placed in a disk of hyperbolic radius $r$.  Arguing
as in Example~\ref{eg:ginibre} we see that $f$ is Lipschitz
with constant $c_r$.  The mean number of points in $\ZZ_\rho$
is $\rho^2 / (1 - \rho^2)$ which for simplicity we can bound
from above by $1/(1-\rho^2)$.  An application of
Theorem~\ref{th:continuous} to $f / c_r$ now yields
$$\P (|f - \E f| \geq a) \leq 5 \exp \left (
   - \frac{a^2}{16 c_r a + 32 c_r^2 (1 - \rho^2)^{-1} } \right ) \, .$$
\end{example}

\section*{Acknowledgments}

We thank Subhroshekhar Ghosh and an anonymous referee for
detailed comments leading to the correction and improvement of
the exposition and results.

\bibliographystyle{alpha}
\bibliography{RP}

\end{document}